\documentclass[10pt]{amsart}

\usepackage{amsmath,amssymb,amsfonts,bbold,mathrsfs,amscd,bm}
\usepackage{pinlabel}

\DeclareMathAlphabet\oldmathcal{OMS}        {cmsy}{b}{n}
\SetMathAlphabet    \oldmathcal{normal}{OMS}{cmsy}{m}{n}
\DeclareMathAlphabet\oldmathbcal{OMS}       {cmsy}{b}{n}

\usepackage{eucal}

\newtheorem{thm}{Theorem}[section]
\newtheorem{prop}[thm]{Proposition}

\newcounter{num}

\newcommand{\Z}{\mathbb{Z}}

\newcommand{\R}{\mathbb{R}}
\newcommand{\C}{\mathbb{C}}

\def\gb{{\mathfrak b}}

\def\gg{{\mathfrak g}}
\def\gh{{\mathfrak h}}

\def\gl{{\mathfrak l}}

\def\gp{{\mathfrak p}}
\def\gq{{\mathfrak q}}

\def\calk{{\mathcal K}}

\def\calw{{\mathcal W}}

\newcommand{\cps}{\mathbb{CP}}
\newcommand{\ol}[1]{\overline{#1}}
\newcommand{\del}{\partial}

\newcommand{\Pic}{\operatorname{Pic}}

\newcommand{\Ric}{\operatorname{Ric}}

\newcommand{\Sl}{\operatorname{\mathfrak{sl}}}

\newcommand{\SU}{\operatorname{SU}}

\title{Calabi-Yau metrics on canonical bundles of flag varieties}
\author{Craig van Coevering}
\address{Department of Mathematics, Bo\u{g}azi\c{c}i University\\
  34342 Bebek, Istanbul, Turkey}
\email{craig.coevering@boun.edu.tr}
\date{July 10, 2018}

\begin{document}

\begin{abstract}
This note gives a simple formula for the unique asymptotically conical Calabi-Yau metrics on the canonical bundle of a flag variety known to exist by the work of 
R. Goto and others.   This is done by generalizing the well known Calabi Ansatz to general K\"{a}hler classes.
We give some examples of explicit families, in particular, a formula for the two dimensional family of asymptotically conical metrics on the canonical bundle of
$F_{1,2}$.
\end{abstract}

\maketitle

\section{Introduction}\label{sec:intro}

Much work has been done on manifolds with a Ricci-flat K\"{a}hler metrics asymptotic to a cone metric at infinity~\cite{vanC10,Got12,ConHei13}.   By a cone we mean a manifold
$\R_{>0} \times M$ with a complex structure $J$ for which a metric $\ol{g}=dr^2 +r^2 g$ is K\"{a}hler, where $g$ is a metric on $M$ and $r$ is the radial coordinate on $\R_{>0}$.
If we add the vertex, such a manifold $Y$ is known to be a normal affine variety, with an isolated singularity.    We say that $N$ with metric $\hat{g}$ is asymptotic to 
to cone metric $\ol{g}$ if there are compact sets $K\subset N$ and $L\subset Y$, and a diffeomorphism $\phi: Y\setminus L\rightarrow N\setminus K$   so that 
on $Y\setminus L$ the metrics satisfy

\begin{equation}\label{eq:conver}
\| \nabla^j \bigl( \phi^*\hat{g} -\ol{g}\bigr) \| =O(r^{-\gamma-j}) 
\end{equation}

for some decay constant $\gamma>0$, where $\nabla$ is the Riemannian connection on $Y$ and the norm is that induced by $\ol{g}$.   
One also requires a similar convergence of the complex structures, or one can 
consider the stronger case in which $\phi$ is as biholomorphism.    Existence and uniqueness results for Ricci-flat(Calabi-Yau) asymptotically conical K\"{a}hler metrics are 
known in the latter case when $\ol{g}$ is Calabi-Yau.  In this case such Calabi-Yau metrics are proved to exist in~\cite{Got12} in each K\"{a}hler class, and are proved to be unique provided 
they satisfy (\ref{eq:conver}) with $\gamma>0$.   Furthermore, the metrics are known to satisfy (\ref{eq:conver}) with $\gamma= 2n$, $\dim_{\C} N =n$, when 
the K\"{a}hler class is a compactly supported class, and $\gamma=2$ otherwise.   

In this note we give a simple formula for these metrics when $N$ is the total space of the canonical bundle $\mathbf{K}_Z$ of a simply connected homogeneous
K\"{a}hler manifold $Z$.   These manifolds are precisely the flag varieties, i.e. $Z=G^{\C}/P$ where  $G^{\C}$ is a complex semisimple Lie group and $P\subset G^{\C}$ 
is a parabolic subgroup.   We suppose that $G\subset G^{\C}$ is a real form, all of which act transitively on $Z$.   
The $G$-invariant K\"{a}hler forms on $Z$, and more generally the G-invariant closed $(1,1)$-forms, can be described explicitly~\cite{AzaBis03}, and they are uniquely determined 
by their class in $H^2(Z)$.   It follows that the invariant K\"{a}hler metric $\omega$ in $\pi c_1(Z)$ is Einstein, since $\Ric_{\omega} \in 2\pi c_1(Z)$ and is $G$-invariant.  
That is
\begin{equation}\label{eq:K-E}
\Ric_{\omega} =2\omega. 
\end{equation}

It follows that \cite{FutOnWang09,Spa11} $\mathbf{K}_Z^\times$, $\mathbf{K}_Z$ minus the zero section, has a Calabi-Yau cone metric.   Note that the Remmert reduction of $\mathbf{K}_Z$ 
gives the affine variety $Y$, mentioned above, and holomorphic map $\varpi: \mathbf{K}_Z \rightarrow Y$.   So $N=\mathbf{K}_Z$ is a crepant resolution of the 
singularity of $Y$ with exceptional divisor $Z$.

We give an explicit Calabi Ansatz formula for the asymptotically conical Calabi-Yau metric on $N$ in each K\"{a}hler class.

\begin{thm}\label{thm:main}
The unique asymptotically conical Ricci flat Kahler metric on $\mathbf{K}_Z$, for each K\"{a}hler class, is given by
\begin{equation}\label{eq:main}
\hat{\omega} = V(\rho) \pi^* \omega +\pi^* \Theta +2\sqrt{-1} V'(\rho)\del r\wedge\ol{\del}r ,
\end{equation}
where $r=h|v|$, $h$ the $G$-invariant Hermitian metric on $\mathbf{K}_Z$, $\rho =r^2$, where $V(\rho)$ is smooth, and $\Theta$ is semi-positive $G$-invariant $(1,1)$-form on $Z$.

This metric is asymptotic to the Calabi-Yau cone metric with decay rate precisely $\gamma =2n$ if $[\hat{\omega}]$ is proportional to $\pi^* c_1(Z)$, and 
precisely $\gamma=2$ otherwise.
\end{thm}

The K\"{a}hler class is class is compact, i.e. contained in cohomology with compact supports $H^2_c (N,\R)$ precisely when it is proportional to $\pi^* c_1(Z)$.
In this case (\ref{eq:main}) is the well known Calabi Ansatz \cite{Cal79}, see also~\cite{Fut11,Cor18}.

\section{Geometry of flag varieties}\label{sec:geo-flag}

We fix some notation; see~\cite{Kna02,VinGorOn94} for background on Lie groups.  
Let $G^{\C}$ be a complex semisimple Lie group, with Lie algebra $\gg^{\C}$.   Let $B\subset G^{\C}$ be a Borel subgroup with Borel subalgebra
$\gb\subset\gg^{\C}$, and Cartan subalgebra $\gh\subset\gb$ and $\Pi$ the simple system of roots associated to the positive system of roots $\Delta^+$ defined by $\gb$. 
Let $G\subset G^{\C}$ be the compact real form, with Lie algebra $\gg$ constricted from the Cartan-Weyl basis of $\gg^{\C}$.
Recall parabolic subgroups of $G^{\C}$ are the subgroups which contain a Borel subgroup.
Let $P\subset G^{\C}$ be a parabolic subgroup, containing $B$.  The Lie algebra $\gp$ of $P$ is well known to be characterized by a subset $S_{\gp}\subset\Pi$.
Then with respect to the choice of Cartan subalgebra $\gh$ and simple system of roots $\Pi$ we have the following root space decompositions 
\begin{gather}\label{eq:parab}
\gb  = \gh \oplus\underset{\alpha\in\Delta^+}{\bigoplus} \gg_{\alpha} \\
\gp  = \gh\oplus\gb\oplus\underset{\alpha\in\langle S_{\gp}\rangle}{\bigoplus} \gg_{-\alpha} 
\end{gather}
And 
\[ \gg =\gp\oplus \gq \]
where 
\begin{equation} \label{eq:nilpo}
 \gq = \underset{\alpha\in\Delta^+ \setminus\langle S_{\gp}\rangle}{\bigoplus} \gg_{-\alpha}. 
\end{equation} 
Further, note that $\ol{\gq} = \underset{\alpha\in\Delta^+ \setminus\langle S_{\gp}\rangle}{\bigoplus} \gg_{\alpha} \subset\gp$ is the nilradical of $\gp$.   And
\begin{equation*}
\gp =\ol{\gq}\oplus \gl, 
\end{equation*}
where 
\[ \gl =\gh \oplus\underset{\alpha\in\langle S_{\gp}\rangle}{\bigoplus}\bigl( \gg_{\alpha} \oplus \gg_{-\alpha}\bigr) \]
is a Levi compliment.  

Note that $G\cap P =L$ is a compact Lie group, whose Lie algebra $\gl_0$ is a real form of $\gl$, and we have $G^{\C}/P =G/L$ .   Furthermore, 
$G\cap B=T$ is a maximal torus of $G$ and a real form of the Cartan subgroup $H$ defined by $\gh$.   

Let $\Pi =\{\alpha_1 ,\ldots, \alpha_\ell \}$ with $\alpha_i \in S_{\gp}$ for $i=k+1, \ldots, \ell$.  And let $\varpi_1, \ldots, \varpi_\ell$ be the corresponding fundamental 
weights, that is the algebraically integral elements of $\gh^*$ with $\langle \varpi_i , \hat{\alpha}_j \rangle =\delta_{ij}$, where $\check{\alpha}_j$ denote the coroots.  

Given a flag variety $Z=G^{\C}/P$, the affine variety $Y$, $\mathbf{K}_Z^\times$ plus the singularity, with its Calabi-Yau cone metric, and its resolution
$N=\mathbf{K}_Z,\ \varpi: N\rightarrow Y$, have a particularly elegant description.   Let $\lambda\in\gh^*$ be dominant algebraically integral, with 
$\langle\lambda \hat{\alpha}_i \rangle >0$ for $i=1,\ldots, k$, and $=0$ for $i=k+1,\ldots,\ell$.   In other words, $\lambda=\sum_{i=1}^k n_i \varpi_i$ with each $n_i >0$.
The character $\chi^{\lambda} :H\rightarrow\C^*$ can be seen to extend 
to $P$.   Denote by $\C(\lambda)$ to be the complex line with the action of $\chi^{\lambda} :P\rightarrow\C^*$.
The associated line bundle $\mathcal{L}_\lambda:= G^{\C}\times_{P} \C(-\lambda)$ is ample.   Let $L(\lambda)$ be the irreducible representation of $G^{\C}$ with 
highest weight $\lambda$.  By the Borel-Weil theorem we have 
\[ L(\lambda)^* =H^0(Z,\mathcal{O}(\mathcal{L}_\lambda)).\]
One can further show that if $v_{\lambda}\in L(\lambda)$ is a highest weight vector, that the orbit $G^{\C}\dot v_{\lambda} \subset L(\lambda)$ gives and embedding
$Z\subset\mathbb{P}(L(\lambda))$.  In other words, $\mathcal{L}_{\lambda}$ is very ample.   

The tangent bundle is $TZ =G^{\C}\times_P \gg/\gp$.  The anticanonical bundle $\mathbf{K}_Z^{-1}=\Lambda^m TZ =G^{\C}\times_P \bigwedge^m (\gg/\gp)$, and is 
easily seen to be $\mathcal{L}_{\delta_P}$,  where $\delta_P :=\sum_{\alpha\in\Delta^+ \setminus\langle S_{\gp}\rangle} \alpha$.  Note that one can show
the character $\chi^{\delta_P}$ extends to $P$ by showing that it is trivial on the commutator subgroup.   We have an embedding
\begin{equation}\label{eq:embed}
Z\subset\mathbb{P}(L(\delta_P))
\end{equation}
so that $\mathbf{K}_Z$ is the restriction of the tautological line bundle on $\mathbb{P}(L(\delta_P))$.  
Each line bundle $\mathcal{L}_\lambda$ has a natural $G$-invariant Hermitian metric, since $\mathcal{L}_\lambda = G\times_L \C(-\lambda)$ and $\chi^{-\lambda}$ restricted to $L$ 
preserves the usual norm on $\C$.   Provide $L(\delta_P)$ with the unique, up to scale, $G$-invariant Hermitian inner product, which provided a Hermitian metric $h$
on the tautological line bundle on $\mathbb{P}(L(\delta_P))$.    This metric restricts to a $G$-invariant metric on $\mathbf{K}_Z$  in (\ref{eq:embed}), since we have 
$G^{\C}\cdot v_{\delta_P} =G \cdot v_{\delta_P}$.   Thus the Fubini-Study metric $\omega_{FS} =\sqrt{-1}\del\ol{\del}\log h|z|^2$, in local holomorphic coordinates,
restricts to the $G$-invariant K\"{a}hler-Einstein metric $\omega$ in (\ref{eq:K-E}).

In~\cite{AzaBis03} it was observed any closed $(1,1)$-form $\omega$ on $Z$ has  a \emph{quasi-potential}, that is, if $\pi: G^{\C}\rightarrow Z$ is the projection then 
$\pi^* \omega =\sqrt{-1}\del\ol{\del} \phi$ for a smooth real valued $G$-invariant function $\phi$.   If $\rho$ is an irreducible representation of $G^{\C}$ with a 
$G$-invariant Hermitian inner product and $v$ is a highest weight vector, then $\phi :=\log \| \rho(g) v \|$
is a $G$-invariant quasi-potential.   That is, there is a closed $(1,1)$-form $\omega_{\phi}$ on $Z$,  so that $\pi^* \omega_{\rho} =\sqrt{-1}\del\ol{\del}\phi$.
In particular, let $\rho^{\varpi_i}$ be the irreducible representation with highest weight $\varpi_i $ for $i=1,\ldots, k$.  Then 
\[ \phi_{\varpi_i} =\log\| \rho^{\varpi_i} v\| \]
is a quasi-potential.   Note that by the $\del\ol{\del}$-lemma, all $G$-invariant closed $(1,1)$-forms are harmonic.

We summarize important properties of flag varieties.  The second result is due to~\cite{AzaBis03}.
\begin{prop}\label{prop:kah-cone}
The Picard group of $Z=G^{\C}/P$ is
\[ \Pic(Z) =\Z\{ \varpi_1,\ldots,\varpi_k \}.\]
The ample, and very ample, line bundles correspond to elements $\sum_{i=1}^k n_i \varpi_i$ with each $n_i >0$.

The closed $G$-invariant $(1,1)$-forms  on $Z$ are $\omega_\phi$ for a quasi-potential
\begin{equation}
\phi =\sum_{i=1}^k c_i \log\| \rho^{\varpi_i} v\|
\end{equation}
with $c_i \in\R$, and $\omega_{\phi}$ is a K\"{a}hler form precisely when $c_i >0$ for $i=1,\ldots, k$.  

Thus the K\"{a}hler cone $\calk_Z$ of $Z$ is identified with the open face of the Weyl chamber $\calw_{\gp}$ of the positive system of roots spanned by
$\{\varpi_1, \ldots, \varpi_k \}$, so
\[ \calk_Z = \R_{>0} \varpi_1 +\cdots +\R_{>0} \varpi_k .\]

The K\"{a}hler cone $\calk_N$ of $N=\mathbf{K}_Z$, is the pull back of $\calk_Z$ under $\pi: \mathbf{K}_Z \rightarrow Z$.  We have the exact sequence
\[ 0\rightarrow H^2_c(N,\R) \rightarrow H^2(N,\R) \rightarrow H^2(M,\R)\rightarrow 0\]
where $M\subset \mathbf{K}_Z^{\times}$ is the $S^1$ subbundle.   And $H^2_c (Z,\R)$ is the line spanned by $\pi^* c_1(\mathbf{K}_X)$, which is identified with
the line spanned by $\delta_{\gp}$ in $\calw_{\gp}$.
\end{prop}

Let $h_{\omega_{\phi}}$ be the Hermitian metric with K\"{a}hler form $\omega_{\phi}$.   For $X,W  \in T^{1,0} Z$ we have 
$h_{\omega_{\phi}}(X,W) =-\sqrt{-1}\omega_{\phi}(X,\ol{W})$.   
One can show that $h_{\omega_{\phi}}$ has the following description~\cite{AzaBis03} in terms of the Cartan-Weyl basis $X_\alpha \in\gg_{\alpha}$ of $\gg$.
The tangent space of $Z$ at $eP$ can be identified with $\gq$ by the infinitesimal acton of $\gq$.   For each root space $\gg_{\beta} \subset\gq$, we have that
$X_{\beta}$ induces a tangent vector at $eP$.   Then on can show that the $X_\beta$ are orthogonal and 
\begin{equation}\label{eq:Her-met}
  h_{\omega_{\phi}}(X_\beta, X_\beta) = \sum_{i=1}^k \frac{c_i}{2} \varpi_i (\check{\beta}) 
\end{equation}

\section{Main theorem}\label{sec:main-thm}

In this section we fix $Z=G^{\C}/P,,\ \dim_{\C} Z=m$, with the $G$-invariant K\"{a}hler-Einstein metric $\omega$ (\ref{eq:K-E}), and provide $\mathbf{K}_Z$ with the $G$-invariant Hermitian 
metric $h$.  Then define a radial function $r$ on the total space $N$ of $\mathbf{K}_Z$ by $r^2 =h(v,v)$ for $v\in \mathbf{K}_Z$.   Then 
$\omega =\sqrt{-1}\del\ol{\del}\log r$.    Let $\Theta$ be a semi-positive $G$-invariant $(1,1)$-form on $Z$.    Set $t=\log r$, so $\omega =\sqrt{-1}\del\ol{\del} t$.
We consider metrics on $\mathbf{K}_Z^\times$ of the form
\[ \hat{\omega}_{\Theta} =\Theta +\sqrt{-1}\del\ol{\del}F(t). \]
In this section we will write $\Theta$ for $\pi^* \Theta$ to simplify notation.  We have 
\begin{equation} \label{eq:Calabi-An}
\hat{\omega} =\Theta +V(t)\omega +\sqrt{-1} V'(t)\del t\wedge\ol{\del}t 
\end{equation}
where we set $V(t)= F'(t)$.   Set $\eta =d^c \log r^2$, where $d^c =\frac{\sqrt{-1}}{2}(\ol{\del} -\del)$.    Then a simple computation gives
\[ \sqrt{-1}\del t \wedge\ol{\del}t =\frac{\sqrt{-1}}{r^2} \del r\wedge\ol{\del}r =\frac{1}{2}\frac{dr}{r}\wedge\eta. \]
And the volume form
\begin{equation}\label{eq:vol-form}
\hat{\omega}_{\Theta} = \frac{m+1}{2} V' \bigl( V\omega +\Theta\bigr)^m \wedge\frac{dr}{r}\wedge\eta. 
\end{equation}

If we rescale the radial function $\tilde{r} =r^{\frac{1}{m+1}}$, then 
\begin{equation}\label{eq:CY-cone} 
\ol{\omega}_{CY} =\frac{1}{2}\sqrt{-1}\del\ol{\del} \tilde{r}^2 =\frac{1}{m+1}\tilde{r}d\tilde{r}\wedge\eta +\frac{\tilde{r}^2}{m+1}\omega 
\end{equation}
is a Ricci flat cone metric on $\mathbf{K}_Z^\times$.   This the K\"{a}hler form of the cone metric $\ol{g}=d\tilde{r}^2 +\tilde{r}^2 g$ over the Sasakian manifold
$M=\{\tilde{r} =1 \} \subset \mathbf{K}_Z^\times$.   This is the $S^1$ subbundle of $\mathbf{K}_Z^\times$ can be given a natural Sasakian structure, which makes
$(M,g)$ is Sasaki-Einstein, by taking the lift of the K\"{a}hler-Einstein metric $\tilde{\omega}=\frac{1}{m+1} \omega$ with Einstein constant $2m+2$ on $Z$.  
More precisely, 
\[ g=\tilde{\eta}\otimes\tilde{\eta} +\pi^* \tilde{g} \]
where $\tilde{g}$ is the Riemannian metric of $\tilde{\omega}$, and $\tilde{\eta} =d^c \log\tilde{r}^2 =\frac{1}{m+1}\eta$ is a contact form with $\frac{1}{2}d\tilde{\eta}=\tilde{\omega}$.
The cone metric $\ol{g}$ is Ricci flat because its Ricci curvature satisfies
\[ \Ric_{\ol{\omega}_{CY}} =\Ric_{\tilde{\omega}} -(2m+2)\tilde{\omega}.\]
See~\cite{FutOnWang09,Spa11} for more details.  

We also have the K\"{a}hler cone metric with respect to the potential $\frac{r^2}{2}$
\[ \ol{\omega} =\frac{1}{2}\sqrt{-1}\del\ol{\del} r^2 = rdr\wedge\eta +r^2 \omega, \]
and a straight forward computation gives 
\begin{equation}\label{eq:vol-ident}
\ol{\omega}^{m+1} =(m+1)^{m+2} r^{2m} \ol{\omega}_{CY}^{m+1}. 
\end{equation}

Let $\mu =\omega^m$ be the unique, up to scale, $G$-invariant volume form on $Z$.   Diagonalize $\Theta$ at a point with respect to $\omega$ at a point,
with eigenvalues $b_i >0,\ i=1,\ldots, m$, then we have
\begin{equation}\label{eq:prod}
\begin{split}
\bigl( V\omega +\Theta\bigr)^m & =\sum_{k=0}^m \binom{m}{k}V^k \omega^k \wedge\Theta \\
	&\sum_{k=0}^m V^k \sigma_{m-k} \mu  \\
\end{split}
\end{equation}
where $\sigma_{m-i}$ is the elementary symmetric function in the $b_i \ i=1,\ldots, m$.   

From (\ref{eq:vol-form}) and (\ref{eq:vol-ident}) we have
\[  \hat{\omega}_{\Theta} = 2^{m-1}(m+1)^{(m+2)} r^{-2}\Bigl( V' \sum_{k=0}^m V^k \sigma_{m-k} \Bigr) \ol{\omega}_{CY}^{m+1} \]
For convenience we will use the variable $\rho =r^2$, so $V'(t)=V'(\rho)2\rho$.
Then $ \hat{\omega}_{\Theta}$ is Ricci flat if 
\begin{equation}\label{eq: Ricci-flat}
V'(\rho)\prod_{j=1}^m (V+b_j ) =C
\end{equation}
for some constant $C>0$.  
Since each $b_i >0$, if we choose $V_0 >0$, then 
\begin{equation}\label{eq:Ricci-flat-int}
 \int_{V_0}^V \prod_{j=1}^m (x+b_j )\, dx =C\rho 
\end{equation}
defines an analytic increasing function $V(\rho)$.  

The equation for (\ref{eq:Ricci-flat-int}) can be written as 
\begin{equation}\label{eq:Ricci-flat-pol}
\sum_{k=0}^m \frac{1}{k+1}V^{k+1} \sigma_{m-k} =C\rho +C_0
\end{equation}

We rewrite (\ref{eq:Calabi-An}) in the variable $\rho =r^2$ to get
\[ \hat{\omega} =\Theta +V(\rho)\omega +\sqrt{-1} 2V'(\rho)\del r\wedge\ol{\del}r \]
which we observe extends to  a non-singular metric on $N=\mathbf{K}_Z$.

Note that when $\Theta=0$ in equation (\ref{eq:Ricci-flat-pol}) we have $\sigma_{m-k}=0$ for $k\neq m$ and
\begin{equation}\label{eq:comp-pot}
V =(m+1)^{\frac{1}{m+1}} \bigl( C\rho +C_0 \bigr)^{\frac{1}{m+1}} 
\end{equation}
which is the familiar Calabi Ansatz for the anti-canonical polarization.

We want to compare $\hat{\omega}$ on the end of $N$ with the Ricci-flat cone metric (\ref{eq:CY-cone}).
We consider the Puiseux series of $V$ in the variable $\tau=\rho^{-1}$, from the polynomial equation
\begin{equation}\label{eq:Ricci-flat-pol2}
\tau\Bigl(\sum_{k=0}^m \frac{1}{k+1}V^{k+1} \sigma_{m-k} -C_0 \Bigr)=C 
\end{equation}
we get the series Puiseux series
\begin{equation}\label{eq:Puiseux}
V(\rho) =\sum_{k=-1}^\infty c_k \rho^{-\frac{k}{m+1}} 
\end{equation}
which converges for large $\rho$, and has first coefficients $c_{-1}=(m+1)^{\frac{1}{m+1}} C^{\frac{1}{m+1}}$ and $c_0 =\frac{-\sigma_1}{m}$.
In the potential $\tilde{r}^2 =r^\frac{2}{m+1}$ we have
\begin{equation}\label{eq:pot-puis}
V(\rho ) =\sum_{k=-1}^\infty c_k \tilde{r}^{-2k} 
\end{equation}
and 
\begin{equation}\label{eq:pot-der-puis}
V'(\rho) = \sum_{k=-1}^\infty c_k \frac{-k}{m+1} \tilde{r}^{-2(m+1+k)}
\end{equation}
Using (\ref{eq:pot-puis}) and (\ref{eq:pot-der-puis}) we have, where $c=(m+1)c_{-1}$,
\[ \begin{split}
\hat{\omega}_{\Theta} -c \ol{\omega}_{CY} & = V(\rho)(m+1)\tilde{\omega} +\Theta +V'(\rho)\tilde{r}^{2m} (m+1)^2 \tilde{r}d\tilde{r}\wedge\tilde{\eta}- c\bigl( \tilde{r}d\tilde{r}\wedge\tilde{\eta}+\tilde{r}^2 \tilde{\omega}\bigr) \\
 &  = (m+1)\bigl( c_0 + O(\tilde{r}^{-2})\bigr)\tilde{\omega} +\Theta + O(\tilde{r}^{-4}) \tilde{r}d\tilde{r}\wedge\tilde{\eta}.
\end{split}\]
Taking the norm of this with the Calabi-Yau cone metric $\ol{g}$ gives (\ref{eq:conver}) with $\gamma =2$.   The asymptotics of the derivatives
in (\ref{eq:conver}) also follow from (\ref{eq:pot-puis}) and (\ref{eq:pot-der-puis}).

Note that the rate $\gamma =2$ cannot be increased unless 
$[\hat{\omega}_{\Theta}] \in H^2_c (Y,\R)$.  Suppose this is the case, so we may take $\Theta =0$ and the potential $V(\rho)$ is given by (\ref{eq:comp-pot}).
The Newton series gives
\[ V(\rho) = (m+1)^{\frac{1}{m+1}}C^{\frac{1}{m+1}} \tilde{r}^2 + (m+1)^{\frac{-m}{m+1}}C^{\frac{-m}{m+1}}C_0 \tilde{r}^{-2m}  +  O\bigl(\tilde{r}^{-4m-2} \bigr) . \]
And also
\[  V'(\rho) = (m+1)^{\frac{-m}{m+1}}C^{\frac{1}{m+1}} \tilde{r}^{-2m} +(-m)(m+1)^{\frac{-2m-1}{m+1}} C^{\frac{-m}{m+1}}C_0 \tilde{r}^{-4m-2} +O\bigl(\tilde{r}^{-6m-4}\bigr). \]
From these formulas we have
\[ \begin{split}
\hat{\omega} -c \ol{\omega}_{CY} & = (m+1)V(\rho)\tilde{\omega} +V'(\rho)\tilde{r}^{2m}(m+1)^2  \tilde{r}d\tilde{r}\wedge\tilde{\eta}- c\bigl( \tilde{r}d\tilde{r}\wedge\tilde{\eta}+\tilde{r}^2 \tilde{\omega}\bigr) \\
	& = \Bigl( (m+1)^{\frac{1}{m+1}}C^{\frac{-m}{m+1}}C_0 \tilde{r}^{-2m} +O\bigl(\tilde{r}^{-4m-2} \bigr)\Bigr)\tilde{\omega}\\
	 &   + \Bigl( (-m)(m+1)^{\frac{-2m-1}{m+1}}C^{\frac{-m}{m+1}}C_0 \tilde{r}^{-2m-2} +O\bigl(\tilde{r}^{-4m-4}  \bigr)\Bigr) \tilde{r}d\tilde{r}\wedge\tilde{\eta}.
\end{split}\]
When we take the norm of this with the Calabi-Yau cone metric $\ol{g}$ we get (\ref{eq:conver}) with $\gamma =2m+2$.

\section{Examples}\label{sec:examp}

We consider examples in which the formula in Theorem~\ref{thm:main} gives information, that is examples of $Z=G^{\C}/P$ in which the K\"{a}hler cone has
dimension greater than one.   We list the three cases in which the formula for $V$ can be solved by radicals.   The theory of flag varieties in Section~\ref{sec:geo-flag} can
be used to express the metrics in coordinates.  A flag variety has an open dense holomorphic coordinate chart given by the action of the nilpotent Lie group $Q$ with Lie algebra 
$\gq$ of (\ref{eq:nilpo}).  If $[e]=e\cdot P \in G^{\C}/P$ then $Q\ni q \mapsto q[e]$ defines a coordinate chart.  

\subsection{Complete flag manifolds}

\vspace{10pt}

Let $Z=G^{\C}/B$, so $S=\emptyset$.   Then $m=\dim_{\C} Z =|\Delta^+|$ and the K\"{a}hler cone $\calk_Z$ is given by the Weyl chamber $\calw$ fixed by the positive system.
The compact K\"{a}hler classes correspond to the ray $\R_{>0} \delta\subset\calw$ spanned by $\delta=\frac{1}{2}\sum_{\alpha\in\Delta^+}\alpha$.   
Note that $\delta$ is integral, $\langle\delta, \check{\alpha}\rangle =1$ for each
$\alpha\in\Pi$, and corresponds to $\mathbf{K}^{-\frac{1}{2}}_Z$, which is very ample by Proposition~\ref{prop:kah-cone}, and $Z$ has Fano index 2.

\subsection{ $\pmb{Z=\cps^1 \times\cps^1}$}\label{ex:1}

A $G$-invariant semi-positive $(1,1)$-form on $Z$ is $\Theta =b_1 \theta_1 +b_2 \theta_2$ where $\theta_1$ and $\theta_2$ are pullbacks of the invariant metrics on the 
$\cps^1$ factors, with radius $\frac{\sqrt{2}}{2}$.   We have $\sigma_1 =b_1 +b_2$ and $\sigma_2 =b_1 b_2$. 

The potential $V(\rho)$ is a solution to 
\[ V^3 +\frac{3}{2}V^2 \sigma_1 +3V\sigma_2 -f(\rho)= 0,\]
where $f(\rho)= C\rho +C_0$.   A formula for a solution to a cubic equation was first published in 1545 by Gerolarmo Cardano~\cite{Car93}.

Define
\begin{align*}
\Delta_0 (b_1, b_2) & = \frac{9}{4}\sigma_1^2 -9\sigma_2 =\frac{9}{4}(b_1 -b_2)^2 \\
\Delta_1 (b_1, b_2, \rho) & = \frac{27}{4}\sigma_1^3 -\frac{81}{2}\sigma_1 \sigma_2 -27f(\rho) \\
	& = \frac{27}{4}\bigl( b_1^3 +b_2^3 -3b_1^2 b_2 -3 b_1 b_2^2 \bigr) -27f(\rho).
\end{align*}
Them
\begin{equation}\label{eq:cube-pot}
V(\rho) =-\frac{1}{2}\sigma_1 +\frac{1}{3(2)^{\frac{1}{3}}}\Bigl( -\Delta_1 +\sqrt{\Delta_1^2 -4\Delta_0^3} \Bigr)^{1/3} + \frac{2^{\frac{1}{3}}}{3}\Delta_0\Bigl( -\Delta_1 +\sqrt{\Delta_1^2 -4\Delta_0^3} \Bigr)^{-1/3} 
\end{equation}
where we choose $C_0 >0$ large enough so that $V(\rho)$ is positive.   Then $V(\rho)$ in (\ref{eq:main}) gives an explicit 2-dimensional family of Calabi-Yau metrics.

The anticanonical embedding of $Z=\cps^1 \times\cps^1 \hookrightarrow\mathbb{P}(W)$ is induced by the orbit of the highest weight vector in the representation $W$
with highest weight $2\delta=2\varpi_1 +2\varpi_2$.   This representation is $W=\Sl(2,\C)\hat{\otimes}\Sl(2,\C)$.  We have the coordinate system $U(z_1 ,z_2)$ on $Z$ given by
the action of 
\[  Q=\Bigl\{\begin{bmatrix} 1 & 0 \\ z_1 &1 \end{bmatrix}\ |\ z_1 \in\C\ \Bigr\} \times \Bigl\{\begin{bmatrix} 1 & 0 \\ z_2 &1 \end{bmatrix}\ |\ z_2 \in\C\ \Bigr\} \]
on the highest weight vector 
\[ Q\cdot \Bigl(\begin{bmatrix} 0 & 1 \\ 0 & 0\end{bmatrix} \otimes \begin{bmatrix} 0 & 1 \\ 0 & 0 \end{bmatrix} \Bigr)= 
\begin{bmatrix} -z_1 & 1 \\ -z_1^2 & z_1 \end{bmatrix} \otimes \begin{bmatrix} -z_2 & 1 \\ -z_2^2  & z_2 \end{bmatrix}. \]
Since the Hermitian metric $h$ on $\mathbf{K}_Z$ is the restriction of that on $W$, if we denote by $w$ the fiber coordinate on $\mathbf{K}_Z$ then over $U$
we have 
\begin{equation}\label{eq:radial2} 
\rho =r^2 =|w|^2 (1 +2|z_1|^2 +|z_1|^4 )(1 +2|z_2|^2 +|z_2|^4 ). 
\end{equation}
The semi-positive $G$ invariant $(1,1)$ forms on $Z$ are 
\begin{equation}\label{eq:Theta2}
\Theta = b_1 \frac{dz_1 \wedge d\ol{z}_1}{(1+|z_1|^2 )^2} + b_1 \frac{dz_2 \wedge d\ol{z}_2}{(1+|z_2|^2 )^2} .
\end{equation}
Combining (\ref{eq:cube-pot}), (\ref{eq:radial2}) and (\ref{eq:Theta2}) in formula (\ref{eq:main}) gives and explicit formula for this metric in $U$.

\subsection{$\pmb{Z=\cps^1 \times\cps^1 \times\cps^1}$}

Similar to the last example the semi-positive $G$-invariant forms are $\Theta =b_1 \theta_1 +b_2 \theta_2 +b_3 \theta_3$, for $b_1,b_2,b_3 \geq 0$.
The potential $V(\rho)$ is a solution to 
\begin{equation}\label{eq:pot-deg4}
 V^4 +\frac{4}{3}V^3 \sigma_1 +2V^2\sigma_2 +4V\sigma_3-f(\rho)= 0,
\end{equation}
where $f(\rho)= C\rho +C_0$.   
We have
\begin{align*}
\sigma_1 & =b_1 +b_2 +b_3  \\
\sigma_2 & = b_1 b_2 + b_2 b_3 + b_1 b_3 \\
\sigma_3 & = b_1 b_2 b_3
\end{align*}
The solution to a quartic equation is credited to Lodovico Ferrari in 1540, but a solution was first published in 1545 by Gerolarmo Cardano~\cite{Car93}.
We get an explicit, though admittedly unwieldy formula in this case.
We first make a substitution converting (\ref{eq:pot-deg4}) to the depressed quartic
\[x^4 +px^2 +qx +r=0,\]
where 
\begin{align*}
p & =  2\sigma_2 -\frac{2}{3}\sigma_1^2  \leq 0\\
q & =  \frac{8}{27}\sigma_1^3 -\frac{4}{3}\sigma_1 \sigma_2 +4\sigma_3   \\
r & =  \frac{-1}{27}\sigma_1^4 -f(\rho) -\frac{4}{3}\sigma_1 \sigma_3 +\frac{2}{9}\sigma_1^2 \sigma_2   \\
\end{align*}
and $V=x-\frac{\sigma_1}{3}$.     We consider the \emph{resolvent cubic}
\[ R(y)=8 y^3 +8py^2 + (2p^2 -8r)y -q^2 =0.\]
Notice that $r \rightarrow -\infty$ as $\rho\rightarrow\infty$.  
If $r^2 \geq -2p$ then one can easily show that $R(\frac{q^2}{-r})>0$.  So by Rolle's theorem 
we have that $R(y)$ has a real root $y_0$ with $0<y_0 < \frac{q^2}{-r}$.   One can solve for $y_0$ in terms of radicals using Cardano's formula
in Example~\ref{ex:1}.   To solve for $y_0$ we define
\begin{align*}
\Delta_0 & =4\sigma_2^2 -16\sigma_1 \sigma_3 -12 f(\rho)  \\
\Delta_1 & =16\sigma_2^3 -96\sigma_1 \sigma_2 \sigma_3 - 12\sigma_1^2 f(\rho) + 432\sigma_3^2 +144\sigma_2 f(\rho) 
\end{align*}
Note that $\Delta_0 <0$ for large $\rho$.   Then
\begin{equation}
y_0 =\frac{-1}{3} +\frac{1}{6 (2)^{\frac{1}{3}}}\Bigl(\Delta_1 +\sqrt{\Delta_1^2 -4\Delta_0^3} \Bigr)^{1/3} +\frac{2^{\frac{1}{3}}}{6}\Delta_0 \Bigl(\Delta_1 +\sqrt{\Delta_1^2 -4\Delta_0^3} \Bigr)^{-1/3}.
\end{equation}
Then the greatest root of (\ref{eq:pot-deg4}) is
\begin{equation}
V(\rho) =\frac{-\sigma_1}{3} - \sqrt{\frac{y_0}{2}} +\frac{1}{2}\sqrt{-2y_0  -2p + q \sqrt{\frac{y_0}{2}}}
\end{equation}
or
\begin{equation}
V(\rho) =\frac{-\sigma_1}{3} + \sqrt{\frac{y_0}{2}} +\frac{1}{2}\sqrt{-2y_0  -2p - q \sqrt{\frac{y_0}{2}}}, 
\end{equation}
depending on whether $q>0$ or $q<0$.

\subsection{$\pmb{Z=F_{1,2} =\SU(3)/ {T^2}}$}
\medskip

This flag manifold is well known as the twistor space of $\ol{\cps}^2$ with the Fubini-Study metric with the opposite orientation, which is anti-self-dual.   
It is one of only two K\"{a}hler twistor spaces of anti-self-dual 4-manifolds, the other being $\cps^3$, the twistor space of $S^2$.
The Ricci-flat cone $Y=C(M)\cup\{o\}$ is actually hyperk\"{a}hler.   So $N=\mathbf{K}_Z$ with the Ricci flat metrics of Theorem~\ref{thm:main} is a resolution of a singular 
hyperk\"{a}hler metric.  

The K\"{a}hler cone of $Z$ is two dimensional, spanned by $\varpi_1 ,\varpi_2$.   The eigenvalues of a $G$-invariant semi-positive form $\Theta$, relative to $\omega$, 
from (\ref{eq:Her-met}) are $b_1 ,b_2, \frac{1}{2}(b_1 +b_2)$.   The potential $V$ is again a solution to (\ref{eq:pot-deg4}) but with 
\begin{align*}
\sigma_1 & = \frac{3}{2}(b_1 +b_2) \\
\sigma_2 & = 2b_1 b_2  +\frac{1}{2}(b_1^2 + b_2^2) \\
\sigma_3 & =\frac{1}{2}(b_1^2 b_2 + b_1 b_2^2)
\end{align*}
And again we make a substitutions converting (\ref{eq:pot-deg4}) to the depressed quartic 
\[x^4 +px^2 +qx +r=0,\] 
where straight forward computation gives
\begin{align*}
p & = 2\sigma_2 -\frac{2}{3}\sigma_1^2 \\
   & = -\frac{1}{2}(b_1 - b_2)^2 \\
q & = \frac{8}{27}\sigma_1^3 -\frac{4}{3}\sigma_1 \sigma_2 +4\sigma_3 \\
   & = 0 \\
r  & = \frac{-1}{27}\sigma_1^4 -\frac{4}{3}\sigma_1 \sigma_3 +\frac{2}{9}\sigma_2 \sigma_1^2 -f(\rho)  \\
   & = \frac{1}{16}(b_1 -b_2 )^4  -b_1^2 b_2^2 -f(\rho)  \\
\end{align*}
The variable are related by $V=x-\frac{\sigma_1}{3}$.   Luckily $q=0$ and the depressed quartic can be easily factored to give 
\[ x=\sqrt{ \frac{-p}{2}+\sqrt{p^2 -4r}}  \]
Substituting the above formulae into this give an expression for the potential $V$ in terms of the parameters $b_1 ,b_2$
\begin{equation}\label{eq:flag-pot}
V(\rho) = \frac{-(b_1 +b_2)}{2} +\sqrt{ \frac{(b_1 -b_2)^2}{4} +2\sqrt{ b_1^2 b_2^2 +C\rho +C_0}}, 
\end{equation}
where $C,C_0 >0$.

We can express the family of Calabi-Yau metrics on $\mathbf{K}_{F_{1,2}}$ explicitly in an open dense holomorphic coordinate chart.
The anticanonical embedding of $Z$ is given as an orbit of $G^{\C}$ in the representation of highest weight $2\varpi_1 +2\varpi_2$.
The fundamental weight $\varpi_1$ is the highest weight of $\Sl(3,\C)$ acting on $\C^3$, while $\varpi_2$ is the highest weight of $\Lambda^2 \C^3$, which
is the dual representation.   By the Weyl dimension formula the representation $W_{m_1, m_2}$ with highest weight $m_1 \varpi_1 +m_2 \varpi_2$ has
\[ \dim_{\C} W_{m_1 ,m_2} =\frac{1}{2}(m_1 +1)(m_2 +1)(m_1 +m_2 +1). \]
The anticanonical embedding $Z\hookrightarrow\mathbb{P}(W_{2,2})$ is the orbit of a highest weight vector by $\Sl(3,\C)$.

We have a coordinate system $U(z_1,z_2,z_3)$ on $Z$ given by the action of 
\[ Q=\Bigl\{\begin{bmatrix} 1 & 0 & 0 \\ z_1 & 1 & 0 \\ z_3 & z_2 & 1 \end{bmatrix} \ |\ z_1 ,z_2 ,z_3 \in\C\ \Bigr\} \]
It suffices to compute the orbit of the highest weight vector of $W_{1,0}^{\otimes^2} \otimes W_{0,1}^{\otimes^2}$, since $W_{2,2}$ is the irreducible component 
containing the highest weight vectors.   If $w$ is the fiber coordinate on $\mathbf{K}_Z$ restricted to $U$, then we compute
\begin{equation}\label{eq:radial4}
\rho =r^2 =|w|^2 \bigl(1+|z_1|^2 +|z_3|^2 \bigr)^2 \bigl(1+|z_2|^2 +|z_1 z_2 -z_3 |^2 \bigr)^2
\end{equation}
The semi-positive $G$-invariant $(1,1)$-forms are expressed in $U$ by 
\begin{equation}\label{eq:Theta4}
\Theta =\sqrt{-1} b_1 \del\ol{\del}\log\bigl(1 +|z_1|^2 +|z_3|^2 \bigr) + \sqrt{-1} b_2 \del\ol{\del}\log\bigl(1 +|z_2|^2 +|z_1 z_2 -z_3|^2 \bigr) 
\end{equation}
Combining (\ref{eq:flag-pot}), (\ref{eq:radial4}) and (\ref{eq:Theta4}) in formula (\ref{eq:main}) gives and explicit formula for this metric in $U$.

\bibliographystyle{amsplain}

\providecommand{\bysame}{\leavevmode\hbox to3em{\hrulefill}\thinspace}
\providecommand{\MR}{\relax\ifhmode\unskip\space\fi MR }
\providecommand{\MRhref}[2]{%
  \href{http://www.ams.org/mathscinet-getitem?mr=#1}{#2}
}
\providecommand{\href}[2]{#2}

\end{document}